\documentclass[12pt]{amsart}
\usepackage{amscd,amssymb,graphics,a4wide,color,hyperref}

\usepackage{mathrsfs}

\newtheorem{theorem}{Theorem}[section]
\newtheorem{lemma}[theorem]{Lemma}
\newtheorem{proposition}[theorem]{Proposition}

\theoremstyle{definition}

\theoremstyle{remark}
\newtheorem{remark}[theorem]{Remark}

\newcommand{\HH} {\mathbb{H}}

\newcommand{\ZZ} {\mathbb{Z}}

\newcommand{\RR} {\mathbb{R}}
\newcommand{\CC} {\mathbb{C}}

\newcommand{\PP} {\mathbb{P}}

\newcommand {\shExt} {\mathcal{E} \!\text{\textit{xt}}}

\newcommand {\shHom} {\mathcal{H}\!\text{\textit{om}}}

\newcommand {\shN}  {\mathcal{N}}

\newcommand {\shT}  {\mathcal{T}}

\newcommand {\dual} {\vee}

\newcommand {\Ext}  {\operatorname{Ext}}

\newcommand {\Hom}  {\operatorname{Hom}}

\renewcommand{\O}  {\mathcal{O}}

\newcommand {\Pic}  {\operatorname{Pic}}

\newcommand {\rank} {\operatorname{rank}}

\newcommand {\Sing} {\operatorname{Sing}}

\newcommand {\T} {\shT}

\def\mydate{\ifcase\month \or January\or February\or March\or
April\or May\or June\or July\or August\or September\or October\or 
November\or December\fi \space\number\day,\space\number\year}

\begin{document}
\def\mapright#1{\smash{
  \mathop{\longrightarrow}\limits^{#1}}}
\def\mapleft#1{\smash{
  \mathop{\longleftarrow}\limits^{#1}}}
\def\exact#1#2#3{0\rightarrow#1\rightarrow#2\rightarrow#3\rightarrow0}
\def\mapup#1{\Big\uparrow
   \rlap{$\vcenter{\hbox{$\scriptstyle#1$}}$}}
\def\mapdown#1{\Big\downarrow
   \rlap{$\vcenter{\hbox{$\scriptstyle#1$}}$}}
\def\dual#1{{#1}^{\scriptscriptstyle \vee}}
\def\invlim{\mathop{\rm lim}\limits_{\longleftarrow}}
\def\rto{\raise.5ex\hbox{$\scriptscriptstyle ---\!\!\!>$}}

\input epsf.tex
\title[A Calabi-Yau threefold with Brauer group $(\ZZ/8\ZZ)^2$]
{A Calabi-Yau threefold with Brauer group $(\ZZ/8\ZZ)^2$}

\author{Mark Gross} \address{Department of Mathematics,
UCSD, 9500 Gilman Drive, La Jolla, CA 92093-0112, USA}
\email{mgross@math.ucsd.edu}
\thanks{This work was partially supported by NSF grant
0204326 and 0505325.}

\author{Simone Pavanelli} 
\address{Nextra Investment Management SGR, Piazza Cadorna 3,
20123 Milano, Italy}
\email{simone\_pavanelli@hotmail.com}

\maketitle
\bigskip
\section*{Introduction}

If $X$ is a non-singular Calabi-Yau threefold with $b_1(X)=0$,
then the analytic
cohomological Brauer group $H^2_{an}(X,\O_X^*)$ is isomorphic
to $H^3(X,\ZZ)_{tors}\cong H^4(X,\ZZ)_{tors}\cong H_2(X,\ZZ)_{tors}$, 
(and this coincides with the cohomological Brauer
group in the \'etale topology). There has been recent interest in finding
examples of Calabi-Yau threefolds with non-trivial Brauer group motivated
partly by string theory (\cite{AspMorr}), partly by their
role in mirror symmetry (\cite{AspMorr},\cite{SlagII}),  and partly by the study
of twisted derived categories (\cite{Cald1}, \cite{Cald2}). 
For a long time, there seemed
to be no examples of non-trivial Brauer groups for simply-connected,
non-singular Calabi-Yau threefolds. Recent work of Batyrev and
Kreuzer \cite{BK} identified all Calabi-Yau hypersurfaces in toric varieties
with non-trivial Brauer group, and examples were found with groups
of orders $2$, $3$, and $5$. In this note we study a Calabi-Yau
threefold first discovered by the first author and Sorin Popescu in 
\cite{GP}, and show that it has Brauer group $(\ZZ/8\ZZ)^2$.

This Calabi-Yau threefold was studied extensively in the second
author's thesis \cite{Pav}, motivated by the hope that it would
be self-mirror. Carrying out the $A$ and $B$ model
instanton calculations provided supporting evidence for this idea;
however, a deeper analysis using the Strominger-Yau-Zaslow conjecture
led to the expectation that the mirror was not simply connected.
These results could be explained if the original Calabi-Yau had 
a non-trivial Brauer group, and is the universal cover of its mirror.
While this has not yet been proven, this proposal led to the calculation
of the Brauer group presented here.

Unfortunately, calculating torsion groups of cohomology is often
difficult, and we apologize if the proof given here is unnecessarily
complicated.

\section{The example}

We recall the construction of a Calabi-Yau threefold $V_{8,y}\subseteq
\PP^7$ from \cite{GP}, \S 6. Here $y$ is a parameter lying in a projective plane
we call $\PP^2_-\subseteq \PP^7$ given by the equations 
\[
x_0=x_1+x_7=x_2+x_6=x_3+x_5=x_4=0.
\]
Let $\sigma,\tau:k[x_0,\ldots,x_7]\rightarrow k[x_0,\ldots,x_7]$ be
the automorphisms given by $\sigma(x_i)=x_{i-1}$ and $\tau(x_i)=\xi^{-i}x_i$, 
with indices taken modulo $8$ and $\xi$ a primitive $8$th root of unity.
These automorphisms generate the Heisenberg group $\HH_8$, a central extension
of $(\ZZ/8\ZZ)^2$ by $\mu_8$, and the action of $\HH_8$ on $k[x_0,\ldots,x_7]$
induces an action of $\HH_8$ on $\PP^7$ which factors through
the quotient $(\ZZ/8\ZZ)^2$. 
Any abelian surface $A$
with a polarization
of type $(1,8)$ can be mapped to $\PP^7$ in such a way that the action
of $\HH_8$ on $\PP^7$ induces translation on $A$ by elements of the kernel
of the polarization. For a review of this standard theory, see \S 1 of
\cite{GPeqs}. 

Take
\begin{eqnarray*}
f_0&=&x_0^2+x_4^2\\
f_1&=&x_1x_7+x_3x_5\\
f_2&=&x_2x_6.
\end{eqnarray*}
We write $y\in\PP^2_-$ as $y=(y_1:y_2:y_3)$, by which we mean the point
$(0:y_1:y_2:y_3:0:-y_3:-y_2:-y_1)$. Given $y$, set
\[
f=y_1y_3f_0-y_2^2f_1+(y_1^2+y_3^2)f_2,
\]
and define $V_{8,y}$ by the equations
\[
f=\sigma(f)=\sigma^2(f)=\sigma^3(f)=0.
\]
Then we have

\begin{theorem}
\label{review}
Let $y\in \PP^2_-$ be general.
\begin{enumerate}
\item
$V_{8,y}$ is a complete intersection, singular only at $64$ points,
each singularity being an ordinary double point. These $64$ points
are the $\HH_8$-orbit of $y$.
\item
There is a small resolution $\mu:V^1_{8,y}\rightarrow V_{8,y}$ and
a fibration $\pi:V^1_{8,y}\rightarrow\PP^1$ whose general fibre
is an abelian surface with a $(1,8)$ polarization induced by its embedding
in $\PP^7$. (In \cite{GP}, this map is called $\pi_1$.)
\item
$\chi(V^1_{8,y})=0$ and $h^{1,1}(V^1_{8,y})=h^{1,2}(V^1_{8,y})=2$.
\item
$H^2(V^1_{8,y},\ZZ)/Tors$ is generated by $H$, the pull-back of the
class of a hyperplane section of $V_{8,y}$, and $A$ the class
of a fibre of $\pi$. Also $H^4(V^1_{8,y},\ZZ)/Tors$ is generated by
$[l]$, the class of a line in $V_{8,y}$ disjoint from $\Sing(V_{8,y})$,
and $[e]$ the class of an exceptional curve of $\mu$. We have
$A.e=1$, so any exceptional curve of $\mu$ is a section of $\pi$.
In addition,
\[
H^3=16, H^2A=16, A^2=0.
\]
\item
$V_{8,y}\cap \PP^2_-=\{y,\sigma^4(y),\tau^4(y),\sigma^4\tau^4(y)\}$,
and $V_{8,y}$ contains every $\HH_8$-invariant abelian surface $A\subseteq
\PP^7$ with $y\in A$ which is 
embedded with a polarization of type $(1,8)$.
The proper transform of each such abelian surface in $V^1_{8,y}$
is a fibre of $\pi$.
\item
Any $\HH_8$-invariant abelian surface $A\subseteq\PP^7$ embedded with a
polarization of type $(1,8)$ intersects $\PP^2_-$ in four points.

\end{enumerate}
\end{theorem}

\proof This is Theorems 6.5, 6.7, 6.9 and Remark 6.10 of \cite{GP}.
\qed

\medskip

\begin{lemma}
\label{groupscheme}
Let $V^{1,\#}_{8,y}=V^1_{8,y}\setminus Crit(\pi)$, where $Crit(\pi)$
denotes the closed subset where $\pi$ is not smooth. Then 
after fixing any one of the $64$ exceptional curves of $\mu$ to be
the identity section of $\pi$,
$\pi^{\#}:V^{1,\#}_{8,y}
\rightarrow B:=\PP^1$ is a group scheme.
\end{lemma}

\proof 
We will in fact show that $\pi^{\#}$ is a N\'eron
model, hence a group scheme.

We first note that no fibre of $\pi$ can be reducible
or non-reduced: by Theorem \ref{review}, (4), any surface $S\subseteq
V^1_{8,y}$ has $\deg \mu(S)$ divisible by $16$. Since $\deg\mu(A)=16$,
$A$ can't be written as the sum of other effective divisors. Thus each
singular fibre is at worst singular along a curve. 

Working locally at a point $b\in B$, we note that given a 
local section $\sigma$ of $\pi$ over an \'etale neighbourhood $U$ of $b$, 
translation by this section 
induces a birational map $T_{\sigma}:V_U\rto V_U$, where
$V_U=V^1_{8,y}\times_B U$. However, $\pi_U:V_U\rightarrow U$ is relatively
minimal, and hence $T_{\sigma}$ becomes an isomorphism after performing
a sequence of flops on $V_U$ over $U$. On the other hand, the general
fibre of $\pi$ has Picard number 1 and every fibre is irreducible, 
from which it follows that the
relative $N^1(V_U/U)$ is one-dimensional. Here $N^1(V_U/U)=(\Pic V_U\otimes
\RR)/\equiv$, where $D\equiv 0$ if $D.C=0$ for all curves $C$
contracted by $\pi$. This implies that if any curve in a fibre of
$\pi_U$ is contracted, all curves are, and there are no curves in fibres
of $\pi_U$ which can be flopped. Thus $T_{\sigma}$ is biregular.

Given this, one can now adapt the standard proof that a relatively minimal
elliptic fibration over a curve is a N\'eron model after removing non-smooth
points;
see for example \cite{Bosch}, pg. 21,
Prop. 1, or \cite{Artin}, Prop. 1.15, page 218. 
\qed

\medskip

Let $E\subseteq\PP^7$ be an elliptic normal curve (of degree 8).
Then a \emph{translation scroll} over $E$ is obtained by choosing
a point $\rho\in E$, and taking
\[
S_{E,\rho}:=\bigcup_{P\in E} \langle P,P+\rho\rangle,
\]
where $\langle P,P+\rho\rangle$ denotes the line
joining $P$ and $P+\rho$. For general choice of $\rho$,
the translation scroll is singular only
along $E$.

\begin{proposition}
For general choice of $y$, $\pi$ has eight singular fibres, and they
are all translation scrolls.
\end{proposition}

\proof 
Consider 
\[
0\rightarrow \pi^*\Omega^1_{\PP^1}\rightarrow\Omega^1_{V^1_{8,y}}
\rightarrow\Omega^1_{V^1_{8,y}/\PP^1}\rightarrow 0.
\]
Dualizing this gives the exact sequence
\begin{equation}
\label{exact}
0\mapright{} \dual{(\Omega^1_{V^1_{8,y}/\PP^1})}
\mapright{\alpha}\shT_{V^1_{8,y}}\mapright{\pi_*}\pi^*\shT_{\PP^1}.
\end{equation}
Clearly $\pi_*$ drops rank precisely on $Crit(\pi)$. 
Choosing an exceptional curve $e$ as the zero-section of 
$\pi$ as in Lemma \ref{groupscheme},
identifying $e$ with $\PP^1$ via $\pi$, the fact that $\pi^{\#}$ is
a group scheme shows that
the relative tangent bundle $\T_{V^{1,\#}_{8,y}/\PP^1}$ is isomorphic
to $\pi^{\# *}\shN_{e/V^1_{8,y}}$. But $\T_{V^{1,\#}_{8,y}/\PP^1}
=\dual{(\Omega^1_{V^1_{8,y}/\PP^1})}|_{V^{1,\#}_{8,y}}$ and
$\dual{(\Omega^1_{V^1_{8,y}/\PP^1})}$ is reflexive, and hence is equal
to $j_*(\dual{(\Omega^1_{V^1_{8,y}/\PP^1})}|_{V^{1,\#}_{8,y}})$,
where $j:V^{1,\#}_{8,y}\rightarrow V^1_{8,y}$ is the inclusion.
Thus the exact sequence (\ref{exact}) becomes
\[
0\mapright{} \pi^*(\O_{\PP^1}(-1)\oplus\O_{\PP^1}(-1))
\mapright{\alpha}\T_{V^1_{8,y}}\mapright{\pi_*}\pi^*\O_{\PP^1}(2).
\]
On the other hand, $\alpha$ must drop rank precisely when
$\pi_*$ drops rank. Since this locus is at most dimension one,
it follows this happens on a pure codimension one locus, and
the class of this degeneracy locus is $c_2(\T_{V^1_{8,y}}\otimes
\pi^*(\O_{\PP^1}(1)))=c_2(V^1_{8,y})$.
On the other
hand, $\mu_*c_2(V^1_{8,y})$ has degree equal to the degree of 
the second Chern class of a general $(2,2,2,2)$ complete intersection
in $\PP^7$. This degree is $64$, so we find $\deg\mu(Crit(\pi))\le 
64$ (with equality if the correct scheme structure is put on $Crit(\pi)$).

We now show for general $y\in\PP^2_-$, there are at least $8$ $\HH_8$-invariant
translation scrolls containing $y$. By \cite{GPeqs}, Theorem 3.1, any
such translation scroll is a degeneration of $\HH_8$-invariant abelian
surfaces, and hence by Theorem \ref{review}, (5) and (6), 
is a fibre of $\pi$. Since each translation scroll is singular along
an elliptic curve of degree $8$,
this will demonstrate that $V^1_{8,y}$ has precisely $8$ singular fibres,
each a translation scroll.

Consider the $4\times 4$ Moore matrix
\[
M_4=(x_{i+j}y_{i-j}+x_{i+j+4}y_{i-j+4})_{0\le i,j\le 3}.
\]
Here, all indices are modulo $8$.
Then by \cite{GPeqs}, Cor. 2.2, if $E\subseteq\PP^7$ is an $\HH_8$-invariant
elliptic curve, then $M_4$ is rank one on $E\times E\subseteq
\PP^7\times\PP^7$. This shows in particular
that if $S(8)$ denotes the surface which is the closure of the union of all 
$\HH_8$-invariant elliptic normal curves of degree $8$ in $\PP^7$, then 
the rank of the matrix 
\[
(y_{i+j}y_{i-j}+y_{i+j+4}y_{i-j+4})_{0\le i,j\le 3}
\]
is one on $S(8)$. The compactified moduli space of $\HH_8$-invariant
elliptic curves is the compactified moduli space of elliptic curves with
a full level structure of degree $8$. (A full level structure
is an identification of the
$8$-torsion points on $E$ with $(\ZZ/8\ZZ)^2$). This space is the modular
curve $X(8)$; by \cite{Shimura}, \S 1.6, the genus of $X(8)$ is $5$. (We
remark that $S(8)$ is not the universal family over $X(8)$ because it does
not have a well-defined section; rather it is a twist of the universal
family.)

Next, given a general
$E\subseteq\PP^7$ an $\HH_8$-invariant elliptic normal
curve, the union of translation scrolls of $E$ is the secant variety $Sec(E)$,
whose ideal is generated by the $3\times 3$ minors of $M_4$ whenever $y$
is chosen to be a general point of $E$. This is \cite{GPeqs}, Theorem 5.2.
Let us consider the intersection of $Sec(E)$ with $\PP^2_-$: this is
obtained by setting $x_i=-x_{8-i}$ in the Moore matrix,
and $Sec(E)\cap \PP^2_-$ is then the locus where the $3\times 3$ minors
of this latter matrix vanishes. This matrix is
\[
\begin{pmatrix}
0&-x_3y_3+x_1y_7&-x_2y_2+x_2y_6&-x_1y_1+x_3y_5\\
x_1y_1-x_3y_5&x_2y_0-x_2y_4&-x_1y_3+x_3y_7&0\\
x_2y_2-x_2y_6&x_3y_1-x_1y_5&0&x_1y_3-x_3y_7\\
x_3y_3-x_1y_7&0&-x_3y_1+x_1y_5&-x_2y_0+x_2y_4
\end{pmatrix}.
\]
We note that if we interchange the second and fourth rows
of this matrix, it becomes skew-symmetric, and hence the
$3\times 3$ minors vanish exactly where the Pfaffian of this
matrix vanishes. The Pfaffian can be computed as 
\[
w_0(y_1^2-y_3^2+y_5^2-y_7^2)/2+w_1(y_0-y_4)(y_2-y_6)+w_2(y_3y_7-y_1y_5)
\]
where
\[
w_0=2x_1x_3,\quad w_1=-x_2^2,\quad w_2=x_1^2+x_3^2.
\] 
Thus in particular, $Sec(E)\cap\PP^2_-$ is a conic (as $Sec(E)$ cannot contain
all of $\PP^2_-$), and this gives a map $\psi'$ from
$X(8)$ to the space of conics in $\PP^2_-$: for $p\in X(8)$, $E_p$ the
corresponding $\HH_8$-invariant elliptic curve in $\PP^7$, 
$\psi'(p)=Sec(E_p)\cap\PP^2_-$. More generally, the above equations
show how
a general point $y=(y_0:\cdots:y_7)\in\PP^7$ determines a conic
in $\PP^2_-$, and hence
define a rational map $\psi:\PP^7\rto\PP(H^0(\PP^2_-,
\O_{\PP^2_0}(2)))$ to the space of conics in $\PP^2_-$. The image
of this map lies in the $\PP^2$ spanned by the conics
$w_0$, $w_1$ and $w_2$.
Let $\hat w_0,\hat w_1,\hat w_2$ be homogeneous coordinates on this $\PP^2$,
so $\psi^*(\hat w_0)=(y_1^2-y_3^2+y_5^2-y_7^2)/2$ and so on. 
A simple calculation performed on Macaulay or Macaulay 2
(\cite{Mac} or \cite{Mac2}) shows that $\psi^*(\hat w_1^4-8\hat w_0^3\hat w_2-8\hat w_0
\hat w_2^3)$ is contained in the ideal generated by the $2\times 2$ minors
of $M_4(y,y)$. Thus in particular, $\psi(S(8))$ is contained in the
curve $C$ defined by $\hat w_1^4-8\hat w_0^3\hat w_2-8\hat w_0
\hat w_2^3=0$. This shows that $\psi'$ factors through $C$.
The map $\psi'$
is non-constant, as was observed in the proof of \cite{GP}, Theorem 6.3.
As the genus of $C$ is $3$ and the genus of $X(8)$ is $5$, $\psi'$ must
be a double cover. 

Thus, for $y\in \PP^2_-$ general, there are precisely $8$ points
$p_1,\ldots,p_8$ of $X(8)$
such that $y\in \psi'(p_i)$ for each $i$. For each such point $p_i$,
there is at least one $\HH_8$-invariant
translation scroll $S_{E_{p_i},\rho_i}\subseteq
Sec(E_{p_i})$ containing $y$.
This yields the desired $8$ translation scrolls containing $y$.
\qed

\bigskip

We now calculate the Brauer group of $V^1_{8,y}$. 
We note that $\HH_8$ acts on $V_{8,y}$ via the
action on $\PP^7$.
It is easy to check this action is fixed-point free,
and lifts to $V^1_{8,y}$, as it is given by translation
by $8$-torsion sections on $V^1_{8,y}$. The center of
$\HH_8$ acts trivially on $\PP^7$, hence acts trivially on $V^1_{8,y}$.

\begin{theorem}
$V^1_{8,y}$ is simply connected, $H^4(V^1_{8,y},\ZZ)=
\ZZ^2\oplus G$, where $G\cong\HH_8/Z(\HH_8)\cong(\ZZ/8\ZZ)^2$,
and an inclusion 
$G\rightarrow H^4(V^1_{8,y},\ZZ)$ is given by $\varphi\in G
\mapsto [\varphi(e_0)]-[e_0]$, where $e_0$ is any one of the exceptional
curves of $\mu$. Here $\varphi\in G$ acts via the action of $\HH_8$ on
$V^1_{8,y}$.
\end{theorem}

Proof. Let $Y$ be a non-singular $(2,2,2,2)$
complete intersection. Topologically, $Y$ is obtained from $V^1_{8,y}$ by
replacing a neighbourhood diffeomorphic to $B^4\times S^2$ of each exceptional $S^2=\PP^1$ 
with a neighbourhood diffeomorphic to $B^3\times
S^3$. Here $B^n$ denotes an $n$-ball. In particular, each exceptional
$S^2$ in $V^1_{8,y}$ is replaced by an $S^3$ in $Y$. Thus if $E$ is the union
of exceptional $S^2$'s in $V^1_{8,y}$, $F\subseteq Y$ the union of 
$S^3$'s, then $V^1_{8,y}\setminus E$ is homeomorphic to
$Y\setminus F$, and it follows from
Van Kampen's theorem that $\pi_1(V^1_{8,y})\cong 
\pi_1(V^1_{8,y}\setminus E)\cong \pi_1(Y)=0$.

In addition, the relative cohomology exact sequence yields
sequences
\begin{eqnarray*}
0&\rightarrow&H^3(V^1_{8,y},\ZZ)\rightarrow H^3(V^1_{8,y}\setminus E,\ZZ)
\rightarrow H^0(E,\ZZ)\\
&\rightarrow&H^4(V^1_{8,y},\ZZ)\rightarrow H^4(V^1_{8,y}\setminus E,\ZZ)\rightarrow
H^1(E,\ZZ)=0
\end{eqnarray*}
and
\[
0=H^1(F,\ZZ)\rightarrow H^4(Y,\ZZ)\rightarrow H^4(Y\setminus F,\ZZ)
\rightarrow H^2(F,\ZZ)=0.
\]
Thus $\ZZ\cong H^4(Y,\ZZ)\cong H^4(V^1_{8,y},\ZZ)/L$,
where $L$ is the subgroup of $H^4(V^1_{8,y},\ZZ)$ generated by Poincar\'e
dual classes of exceptional curves.

Now consider the Leray spectral sequence for $\pi:V^1_{8,y}\rightarrow B$,
$B=\PP^1$. We have some work to do to get control over the various terms
in this sequence.

We have $R^p\pi_*\ZZ=0$ for $p>4$, $\pi_*\ZZ=\ZZ$ and $R^4\pi_*\ZZ=\ZZ$,
as all fibres are irreducible and have a well-defined orientation.
Note also that if $B_0=B\setminus \{p\in B|\hbox{$\pi^{-1}(p)$ is
singular}\}$, and $i:B_0\hookrightarrow B$ is the inclusion, 
$\pi_0:\pi|_{\pi^{-1}(B_0)}:\pi^{-1}(B_0)\rightarrow B_0$, then 
\[
i_*R^p
\pi_{0*}\ZZ=R^p\pi_*\ZZ
\]
for $0\le p\le 4$. This is the notion of a \emph{simple} torus fibration
introduced in \cite{SlagI}. To see this, note that from the standard
construction of a smoothing of a translation scroll (see for example
\cite{AN} or \cite{HW}), a neighbourhood of the translation
scroll in $X$ is topologically equivalent to $T^2\times X$, where $X\rightarrow
D$ is a $T^2$-fibration over a disk with one degenerate fibre, a pinched
torus (i.e. a Kodaira type $I_1$ fibre). Then in a suitable basis, the
monodromy around this singular fibre is
\[
\begin{pmatrix}
1&1&0&0\\
0&1&0&0\\
0&0&1&0\\
0&0&0&1
\end{pmatrix}
\]
The statement about simplicity then follows as in \cite{TMS}, Theorem 1.3.
In fact, it similarly follows that
\[
i_*R^p
\pi_{0*}\ZZ/n\ZZ=R^p\pi_*\ZZ/n\ZZ
\]
for $0\le p\le 4$ and integer $n$.
We can now argue similarly to the proof of \cite{SlagII}, Theorem 3.9, 
as follows.
Poincar\'e-Verdier duality on the two-dimensional space $B$ tells us that
\[
\RR\Hom(\RR\Gamma(B,R^i\pi_*\ZZ),\ZZ)\cong
\RR\Gamma\RR\shHom(R^i\pi_*\ZZ,\ZZ[2]).
\]
Applying $H^{-j}$ to the complexes on the left
and right-hand side, we obtain
\[
H^{-j}(\RR\Hom(\RR\Gamma(B,R^i\pi_*\ZZ),\ZZ))\cong
\Ext^{2-j}(R^i\pi_*\ZZ,\ZZ).
\]
To compute the left-hand side, use the spectral sequence with
\[
E^{pq}_2=\Ext^p(H^{-q}(B,R^i\pi_*\ZZ),\ZZ)\Rightarrow 
H^n(\RR\Hom(\RR\Gamma(B,R^i\pi_*\ZZ),\ZZ)),
\]
which yields exact sequences
\begin{eqnarray*}
0\rightarrow\Ext^1(H^{j+1}(B,R^i\pi_*\ZZ),\ZZ)
&\rightarrow& H^{-j}(\RR\Hom(\RR\Gamma(B,R^i\pi_*\ZZ),\ZZ))\\
&\rightarrow& \Hom(H^j(B,R^i\pi_*\ZZ),\ZZ)\rightarrow 0.
\end{eqnarray*}
To compute the right-hand side,
we use the local-global $\Ext$ spectral sequence. Note that
fibre-wise Poincar\'e duality for the fibres of $\pi$ over $B_0$
yields 
\[
R^p\pi_*\ZZ\cong\shHom(R^{4-p}\pi_*\ZZ,\ZZ)
\]
on $B_0$. The fact that $i_*R^p\pi_{0*}\ZZ=R^p\pi_*\ZZ$ on $B$
then tells us the same isomorphism holds on $B$. The higher sheaf
$\shExt$'s can be computed locally on a contractible neighbourhood $U$ of
a singular point $p$: let $U_0=U\setminus\{p\}$ and $i:U_0\rightarrow U$
be the inclusion. We compute $\Ext^p_U(R^1\pi_*\ZZ,\ZZ)$ using the
description of monodromy; the other cases are similar. Because there are
locally three independent monodromy invariant sections, we have an exact
sequence on $U$
\[
0\rightarrow \ZZ^3\rightarrow R^1\pi_*\ZZ\rightarrow i_!\ZZ\rightarrow 0.
\]
Applying $\Hom_U(\cdot,\ZZ)$ and the fact that $\Ext^p_U(i_!\ZZ,\ZZ)
=\Ext^p_{U_0}(\ZZ,\ZZ)=H^p(U_0,\ZZ)$, we obtain $\Ext^p_U(R^1\pi_*\ZZ,\ZZ)
=0$ for $p\ge 2$ and the exact sequence
\[
0\rightarrow \ZZ\rightarrow \Hom_U(R^1\pi_*\ZZ,\ZZ)
\rightarrow\ZZ^3\rightarrow H^1(U_0,\ZZ)\rightarrow\Ext^1_U(R^1\pi_*\ZZ,\ZZ)
\rightarrow 0.
\]
From the description of monodromy, it is then clear the map $\ZZ^3
\rightarrow H^1(U_0,\ZZ)=\ZZ$ is surjective and $\Ext^1_U(R^1\pi_*\ZZ,\ZZ)=0$.
Similarly, one finds $\shExt^p_B(R^q\pi_*\ZZ,\ZZ)=0$ for $p>0$ and all $q$.
Thus $\Ext^{2-p}(R^q\pi_*\ZZ,\ZZ)=H^{2-p}(B,\shHom(R^q\pi_*\ZZ,\ZZ))
=H^{2-p}(B,R^{4-q}\pi_*\ZZ)$.

Putting this all together, we see that if
\[
e^{p,q}=\rank_{\ZZ} H^p(B,R^q\pi_*\ZZ),
\]
then 
\[
e^{p,q}=e^{2-p,4-q},
\]
and that if 
\[
T^{p,q}=Tors(H^p(B,R^q\pi_*\ZZ)),
\]
then 
\[
T^{p,q}\cong T^{3-p,4-q}.
\]
Of course, $R^q\pi_*\ZZ$ is torsion free, so 
\[
H^0(B,R^q\pi_*\ZZ)_{tors}=0
\]
anyway.

Next consider $H^0(B,R^3\pi_*\ZZ)=H^0(B_0,R^3\pi_{0*}\ZZ)$. By 
\cite{Deligne}, Cor. 4.1.2, $H^0(B_0,R^3\pi_*\ZZ)$
is a sub-Hodge structure of $H^3(\pi^{-1}(b),\ZZ)$ for any $b\in B_0$,
and the Hodge structure on $H^0(B_0,R^3\pi_*\ZZ)$ is independent of
$b$.
However, if this group is non-zero, this implies $\pi^{-1}(b)$ contains a 
fixed non-trivial abelian subvariety
for any $b\in B_0$. However, since the abelian
surfaces appearing in $\pi:V^1_{8,y}\rightarrow B$ are general, this can't
happen. Thus $H^0(B,R^3\pi_*\ZZ)=0$. The entries in the $E^2$-term
of the Leray spectral sequence of $\pi$ are now
\[
\begin{matrix}
\ZZ&0&\ZZ\\
0&\ZZ^{e^{1,3}}\oplus T^{1,3}&\ZZ^{e^{2,3}}\oplus T^{2,3}\\
\ZZ^{e^{0,2}}&\ZZ^{e^{1,2}}\oplus T^{1,2}&\ZZ^{e^{2,2}}\oplus T^{2,2}\\
\ZZ^{e^{0,1}}&\ZZ^{e^{1,1}}\oplus T^{1,1}&T^{2,1}\\
\ZZ&0&\ZZ
\end{matrix}
\]
We note that because $\pi$ has a section, the boundary map
$d:H^0(B,R^4\pi_*\ZZ)\rightarrow H^2(B,R^3\pi_*\ZZ)$ is zero, so
$H^2(B,R^3\pi_*\ZZ)$ injects into $H^5(V^1_{8,y},\ZZ)=0$, so
$H^2(B,R^3\pi_*\ZZ)=0$. Thus $e^{2,3}=e^{0,1}=0$, and $T^{2,3}=T^{1,1}=0$.
Next, since $b_2(V^1_{8,y})=2$ and $e^{0,2}\ge 1$, with $H$ generating
a non-trivial class in $H^0(B,R^2\pi_*\ZZ)$, we see $e^{1,1}=0$ and
$e^{0,2}=1$. So we now have
\[
\begin{matrix}
\ZZ&0&\ZZ\\
0&T^{1,3}&0\\
\ZZ&\ZZ^{e^{1,2}}\oplus T^{1,2}&\ZZ\oplus T^{2,2}\\
0&0&T^{1,3}\\
\ZZ&0&\ZZ
\end{matrix}
\]
Now $\pi:V^1_{8,y}\rightarrow B$ has a group of sections, once one of the
64 exceptional curves is chosen as a zero-section. This group of
sections contains the group $(\ZZ/8\ZZ)^2$, given by those 64
sections. Note that if $\pi^{-1}(b)=V/\Lambda$, where $V=\CC^2$ and
$\Lambda$ is a four-dimensional lattice, then $H^3(\pi^{-1}(b),\ZZ)=\Lambda$
and $H^3(\pi^{-1}(b),\ZZ/8\ZZ)=\Lambda/8\Lambda\cong {1\over 8}\Lambda/\Lambda$.
Thus $H^3(\pi^{-1}(b),\ZZ/8\ZZ)$ is canonically isomorphic to the set
of 8-torsion points on $\pi^{-1}(b)$, and $H^0(B,R^3\pi_*\ZZ/8\ZZ)
=H^0(B_0,R^3\pi_{0*}\ZZ/8\ZZ)$ is equal to the space of
$8$-torsion sections, and hence contains the subgroup
$(\ZZ/8\ZZ)^2$. However, from the exact sequence
\[
0\mapright{} R^3\pi_*\ZZ \mapright{\cdot 8}
R^3\pi_*\ZZ\mapright{} R^3\pi_*(\ZZ/8\ZZ)\mapright{} 0
\]
we see that $H^0(B,R^3\pi_*\ZZ/8\ZZ)=H^1(B,R^3\pi_*\ZZ)_8$,
where the subscript denotes the part of the group killed by $8$.
Thus $(\ZZ/8\ZZ)^2\subseteq T^{1,3}$. This subgroup is naturally 
identified with $G$.

The Leray spectral sequence now gives us a filtration
\[
0=F^0\subseteq F^1\subseteq F^2\subseteq F^3=H^4(V^1_{8,y},\ZZ)
\]
with $F^1=\ZZ\oplus T^{2,2}$ and $F^2/F^1=T^{1,3}$, so $H^4(V^1_{8,y},\ZZ)$
contains the subgroup $F^2$ which is an extension
\[
0\rightarrow \ZZ\oplus T^{2,2}\rightarrow F^2\rightarrow T^{1,3}\rightarrow 0.
\]
We will show there is a map $\varphi:G\rightarrow F^2$ whose composition
with the projection to $T^{1,3}$ is injective.

We will first have to show $T^{1,2}$ contains no two-torsion. Given $b\in B_0$,
$H^3(\pi^{-1}(b),\ZZ/2\ZZ)=\Lambda/2\Lambda$, and we have the wedge product
\[
(\Lambda/2\Lambda)\times (\Lambda/2\Lambda)\rightarrow
{\bigwedge}^2(\Lambda/2\Lambda)
\]
and
\[
(\Lambda/2\Lambda)\times \left({\bigwedge}^2(\Lambda/2\Lambda)\right)
\rightarrow{\bigwedge}^3(\Lambda/2\Lambda).
\]
Given linearly independent
monodromy invariant elements $e_1,e_2\in\Lambda/2\Lambda$,
we obtain a monodromy invariant element $e_1\wedge e_2$
of $\bigwedge^2(\Lambda/2\Lambda)$. Suppose there is a two-dimensional space
of monodromy elements of $\bigwedge^2(\Lambda/2\Lambda)$,
say generated by $e_1\wedge e_2$ and $f$. Then
either $e_1\wedge f\not=0$ or $e_2\wedge f\not=0$
in $\bigwedge^3(\Lambda/2\Lambda)$, and this would
give a monodromy invariant element of $\bigwedge^3\Lambda/2\bigwedge^3\Lambda$.
Thus we see that the fact that $H^0(B,R^3\pi_*(\ZZ/2\ZZ))$
contains two linearly independent elements implies that if
$\dim_{\ZZ/2\ZZ} H^0(B,R^2\pi_*(\ZZ/2\ZZ))\ge 2$, then
$H^0(B,R^1\pi_*(\ZZ/2\ZZ))\not=0$, a contradiction. Thus
$H^1(B,R^2\pi_*\ZZ)$ contains no $2$-torsion.

Now consider $\alpha\in G$, acting on $V^1_{8,y}$
by translation by an $8$-torsion section. Then $\alpha^*:H^4(V^1_{8,y},\ZZ)
\rightarrow H^4(V^1_{8,y},\ZZ)$ induces the identity on rational cohomology. Look at
\[
\log\alpha^* = (\alpha^*-I)-{1\over 2}(\alpha^*-I)^2.
\]
Indeed $\alpha^*-I$ maps $H^p(B,R^q\pi_*\ZZ)$ to $H^{p+1}(B,R^{q-1}\pi_*\ZZ)$
because $\alpha^*-I$ acts trivially on the cohomology of a fibre of $\pi$.
Thus $(\alpha^*-I)^3=0$, and since $T^{2,2}\cong T^{1,2}$,
$H^2(B,R^2\pi_*\ZZ)_{tors}$ has no two-torsion, and so 
${1\over 2}(\alpha^*-I)^2$ is well-defined.

If $e_0$ is the zero section of $\pi$,
then we see the map
$\varphi:G\rightarrow H^4(V^1_{8,y},\ZZ)$ defined by 
$\alpha\mapsto\log\alpha^*([e_0])$ is a homomorphism:
\begin{eqnarray*}
\varphi(\alpha\circ\beta)&=&(\log\alpha^*\circ\beta^*)[e_0]\\
&=&\log\alpha^*([e_0])+\log\beta^*([e_0])\\
&=&\varphi(\alpha)+\varphi(\beta).
\end{eqnarray*}
The image of $\varphi$ is necessarily $G=(\ZZ/8\ZZ)^2$,
identified under projection to $T^{1,3}$
with the space of $8$-torsion sections of $\pi$
contained in $T^{1,3}$. Thus, since
$T^{2,2}$ has no two-torsion, we see
$\varphi(\alpha)=(\alpha^*-I)([e_0])=[e]-[e_0]$,
where $e=\alpha(e_0)\subseteq V^1_{8,y}$. So the subgroup of
$H^4(V^1_{8,y},\ZZ)$ generated by differences of exceptional curves is isomorphic
to $G$. But as observed at the beginning of the proof, $\{[e]-[e_0]\}$
generates the torsion subgroup of $H^4(V^1_{8,y},\ZZ)$, so 
$Tors(H^4(V^1_{8,y},\ZZ))
=(\ZZ/8\ZZ)^2$. \qed

\medskip

\begin{remark} Calculations performed in the second author's thesis
(\cite{Pav}) suggest that the mirror to $V^1_{8,y}$ has fundamental group
$\ZZ/8\ZZ$ and Brauer group $\ZZ/8\ZZ$. This was derived from a conjectural
topological torus fibration on $V^1_{8,y}$ and a calculation of these
groups for the dual fibration. As yet, we do not have a rigorous construction
of a mirror manifold, but one possible guess would be a quotient of
$V^1_{8,y}$ by a subgroup $\ZZ/8\ZZ\subseteq G$.
\end{remark}

\end{document}